\documentclass{amsart}
\usepackage{a4,amsmath,amssymb,amsfonts,times}

\newcommand{\Z}{\mathbb Z}
\newcommand{\R}{\mathbb R}
\newcommand{\GL}{\mathsf{GL}}
\newcommand{\SC}{{\mathbf\Delta}}

\newcommand{\MM}{{\mathcal M}}
\newcommand{\MS}{{\mathcal S}}
\newcommand{\MP}{{\mathcal P}}
\newcommand{\MD}{{\mathcal D}}
\newcommand{\MV}{{\mathcal V}}

\newcommand{\gldz}{\GL_d(\Z)}
\newcommand{\gldr}{\GL_d(\R)}
\newcommand{\sd}{{\mathcal S}^d}
\newcommand{\sdo}{{\mathcal S}^d_{>0}}
\newcommand{\sdgeo}{{\mathcal S}^d_{\geq 0}}

\renewcommand{\vec}[1]{\boldsymbol{#1}}

\DeclareMathOperator{\relint}{relint}
\DeclareMathOperator{\trace}{trace}
\DeclareMathOperator{\conv}{conv}
\DeclareMathOperator{\cone}{cone}
\DeclareMathOperator{\gradient}{grad}
\DeclareMathOperator{\BR}{BR}
\DeclareMathOperator{\Del}{Del}
\DeclareMathOperator{\Min}{Min}
\DeclareMathOperator{\Co}{Co}

\theoremstyle{definition}
\newtheorem{definition}{Definition}[section]
\newtheorem{proposition}[definition]{Proposition}
\newtheorem{theorem}[definition]{Theorem}
\newtheorem{corollary}[definition]{Corollary}
\newtheorem{problem}[definition]{Problem}
\newtheorem{conjecture}[definition]{Conjecture}
\newtheorem{example}[definition]{Example}

\title{Methods in the Local Theory of Packing and Covering Lattices}

\author{A. Sch\"urmann}
\address{Mathematics Department, University of Magdeburg, 39106 Magdeburg, Germany}
\email{achill@math.uni-magdeburg.de}
\author{F. Vallentin}
\address{Einstein Institute of Mathematics, The Hebrew
University of Jerusalem, Israel} 
\email{vallenti@ma.tum.de}
\thanks{The second author was
supported by the Edmund Landau Center for Research in Mathematical Analysis and
Related Areas, sponsored by the Minerva Foundation (Germany).}

\date{September 26, 2008}

\begin{document}

\begin{abstract}
In this paper we are concerned with three lattice problems: the
lattice packing problem, the lattice covering problem and the lattice
packing-covering problem. One way to find optimal lattices for these
problems is to enumerate all finitely many, locally optimal
lattices. For the lattice packing problem there are two classical
algorithms going back to Minkowski and Voronoi. For the covering and
for the packing-covering problem we propose new algorithms.

Here we give a brief survey about these approaches. We report on some
recent computer based computations where we were able to reproduce and
partially extend the known classification of locally optimal
lattices. Furthermore we found new record breaking covering and
packing-covering lattices. We describe several methods with examples
to show that a lattice is a locally optimal solution to one of the three
problems.
\end{abstract}

\maketitle

\section{Introduction}

Classical problems in the geometry of numbers are the determination of
most economical lattice sphere packings and coverings of the Euclidean
$d$-space~$\R^d$. 
A {\em lattice} $L$ is
a full rank, discrete subgroup of $\R^d$.
Thus there exist matrices 
$A\in \gldr$ with $L=A\Z^d$ which we call {\em bases} of $L$.

If $B^d$ denotes the Euclidean unit ball, then the Minkowski sum $L +
\alpha B^d = \{\vec{v} + \alpha\vec{x} : \vec{v} \in L, \vec{x} \in
B^d\}$, $\alpha \in \R_{>0}$, is a lattice packing if the translates
of $\alpha B^d$ have mutually disjoint interiors and a lattice
covering if $\R^d = L + \alpha B^d$. The packing radius $\lambda(L)$
of a lattice $L$ is given by
\[
\lambda(L)
=
\max\{ \lambda : \mbox{$L + \lambda B^d$ is a lattice packing} \},
\]
and the covering radius $\mu(L)$ by
\[
\mu(L)
=
\min \{ \mu : \mbox{$L + \mu B^d$ is a lattice covering} \}.
\]

For a lattice $L$ we define its {\em determinant} $\det(L)=|\det(A)|$, 
which is independent of the chosen basis.
We consider the following three ``quality measures'' of $L$:
\begin{enumerate}
\item the {\em packing density} 
$\delta(L) = \displaystyle\frac{\lambda(L)^d}{\det(L)} \cdot \kappa_d$,
\item the {\em covering density}
$\Theta(L) = \displaystyle\frac{\mu(L)^d}{\det(L)} \cdot \kappa_d$,
\item the {\em packing-covering constant} 
$\gamma(L) = \displaystyle\frac{\mu(L)}{\lambda(L)}.$
\end{enumerate}

Here $\kappa_d={\pi^{d/2}}/{\Gamma(d/2+1)}$ denotes the volume of the unit ball
$B^d$. So the packing density $\delta(L)$ for instance gives the ratio of 
space covered by spheres in the lattice packing $L+\lambda(L)B^d$.
Note that all three quantities are invariant with respect to a scaling
$\alpha L$ of $L$ with $\alpha\not=0$.
For each of the quantities we consider the problem of
finding extremal lattices
attaining a maximum or minimum respectively.

\begin{problem}[Lattice Packing Problem]
For $d\geq 2$, 
determine $\displaystyle\delta_d = \displaystyle \max_L \delta(L)$
and lattices $L$ attaining it.
\end{problem}

\begin{problem}[Lattice Covering Problem]
For $d\geq 2$, 
determine $\Theta_d  =  \displaystyle \min_L \Theta(L)$
and lattices $L$ attaining it.
\end{problem}

\begin{problem}[Lattice Packing-Covering Problem]
For $d\geq 2$, 
determine $\gamma_d = \displaystyle \min_L \gamma(L)$
and lattices $L$ attaining it.
\end{problem}

Note that all three optima are attained. All three problems have in
common that there exist only finitely many local optima for every
$d$ (see Section \ref{sec:lat-quad} for definitions).  
In this article we want to review some of the major tools
available to find such local extrema and to verify their local
optimality. In Section~\ref{sec:results} we briefly summarize known
results of the three problems and in Section~\ref{sec:lat-quad} we
give a short introduction to the connection of lattices and
positive definite quadratic forms, which gives the framework in that
the problems are usually dealt with.  In
Section~\ref{sec:pack-lat} we describe two classical approaches by
Minkowski and Voronoi to enumerate all local optima of the lattice
packing problem.  In Section~\ref{sec:cov-lat} we are concerned with
local optima of the lattice covering and the lattice packing-covering
problem, which can be treated in parallel. Most of these techniques
are described in greater detail in \cite{sv-2004} and \cite{sv-2004b}.

\section{Known Results} 
\label{sec:results}

\subsection{The Lattice Packing Problem}

The lattice packing problem, arising from the study of positive
definite quadratic forms, is the oldest and most popular of the three
problems and has been considered by many authors in the past.  As
shown in Table~\ref{tab:packing}, the solution to the problem was known for
dimension $d\leq 8$ since 1934. For a description of the extremal root
lattices ${\mathsf A}_d,{\mathsf D}_d$ and ${\mathsf E}_d$ and the
history of the problem we refer the interested reader to the book
\cite{cs-1988}.  Recently, Cohn and Kumar \cite{ck-2004} showed that
the Leech lattice $\Lambda$ gives the unique densest lattice packing
in $\R^{24}$.  Furthermore they showed: The density of any sphere
packing (without restriction to lattices) in $\R^{24}$ cannot exceed
the one given by the Leech lattice by a factor of more than $1+1.65
\cdot 10^{-30}$.

\begin{table}[ht]
\begin{center}
\begin{tabular}{|c|c|c|c|}
\hline
$\mathbf d$ & \textbf{lattice} & \textbf{density $\mathbf \delta_d$} & \textbf{author(s)} \\
\hline
 $2$ & ${\mathsf A}_2$    &  $0.9069\ldots$  &  Lagrange \cite{lagrange-1773}, 1773\\
 $3$ & ${\mathsf A}_3={\mathsf D}_3$    &  $0.7404\ldots$  &  Gau{\ss} \cite{gauss-1840}, 1840\\
 $4$ & ${\mathsf D}_4$    &  $0.6168\ldots$  &  Korkine, Zolotareff \cite{kz-1873}, 1873\\
 $5$ & ${\mathsf D}_5$    &  $0.4652\ldots$  &  Korkine, Zolotareff \cite{kz-1877},1877\\
 $6$ & ${\mathsf E}_6$    &  $0.3729\ldots$  &  Blichfeldt \cite{blichfeldt-1934}, 1934\\
 $7$ & ${\mathsf E}_7$    &  $0.2953\ldots$  &  Blichfeldt \cite{blichfeldt-1934}, 1934\\
 $8$ & ${\mathsf E}_8$    &  $0.2536\ldots$  &  Blichfeldt \cite{blichfeldt-1934}, 1934\\[0.3cm]
 $24$ & $\Lambda$    &  $0.0019\ldots$  &  Cohn, Kumar \cite{ck-2004}, 2004\\
\hline
\end{tabular} 
\medskip
\caption{Optimal lattice packings.}
\label{tab:packing}
\end{center}
\end{table}

\subsection{The Lattice Covering Problem}

The lattice covering problem has only been solved up to dimension $5$.  
Recently, we were able to verify the list of known results
computationally. Even more, we found the complete list of $222$ local
covering optima in dimension~$5$.  Table~\ref{table:latcov} invites to a question
formulated by Ryshkov \cite{ryshkov-1967}, who asked for the lowest
dimension in which ${\mathsf A}^*_d$ gives not the thinnest lattice
covering.  In Section~\ref{sec:cov-lat} we describe the method used
to find a lattice $L^{c}_6\subset \R^6$ with $\Theta(L^{c}_6)=
2.4648\ldots < 2.5511\ldots =\Theta({\mathsf A}^*_6)$.  Thus the
answer to Ryshkov's question is $6$.  It remains an open problem to
prove that the lattice $L^{c}_6$ gives the best lattice covering in
dimension $6$.  We do not even know the exact coordinates of the
lattice yet.  Currently, a complete solution in dimension $d\geq 6$
seems out of reach without completely new methods.

\begin{table}[ht]
\begin{center}
\begin{tabular}{|c|c|c|c|}
\hline
$\mathbf d$ & \textbf{lattice} & \textbf{density $\mathbf \Theta_d$} & \textbf{author(s)} \\
\hline
 $2$ & ${\mathsf A}^*_2$    &  $1.2091\ldots$  &  Kershner \cite{kershner-1939}, 1939\\
 $3$ & ${\mathsf A}^*_3$    &  $1.4635\ldots$  &  Bambah \cite{bambah-1954a}, 1954\\
 $4$ & ${\mathsf A}^*_4$    &  $1.7655\ldots$  &  Delone, Ryshkov \cite{dr-1963}, 1963\\
 $5$ & ${\mathsf A}^*_5$    &  $2.1242\ldots$  &  Ryshkov, Baranovskii \cite{rb-1975}, 1975\\
\hline
\end{tabular} 
\medskip
\caption{Optimal lattice coverings.}
\label{table:latcov}
\end{center}
\end{table}

\subsection{The Lattice Packing-Covering Problem}

As in the case of the lattice covering problem the lattice
packing-covering problem has been solved only for dimensions $d\leq 5$.
And, as in the covering case, we recently were able to verify
these results computationally.  Moreover, we found a new best known
lattice $L^{pc}_6$ in dimension $6$, having a slightly lower
packing-covering constant $\gamma(L^{pc}_6)= 1.4110\ldots$ than the
previously best known one $\gamma({\mathsf{E}_6^\ast})=\sqrt{2}$.

One reason for studying the lattice packing-covering problem is the
open question whether there exists a dimension $d$ 
with $\gamma_d \geq 2$. If so, then any
$d$-dimensional lattice packing with spheres would leave 
space large enough for spheres of the same radius.
This would in particular prove that 
densest sphere packings in dimension~$d$  
are non-lattice packings. 
This phenomenon is likely to be true for large dimensions,
but has not been verified for any $d$ so far.

\begin{table}[ht]
\begin{center}
\begin{tabular}{|c|c|c|c|}
\hline
$\mathbf d$ & \textbf{lattice} & \textbf{density $\mathbf \gamma_d$} & \textbf{author(s)} \\
\hline
 $2$ & ${\mathsf A}^*_2$   &  $1.1547\ldots$  &  Ryshkov \cite{ryshkov-1974}, 1974\\
 $3$ & ${\mathsf A}^*_3$   &  $1.2909\ldots$  &  Ryshkov \cite{ryshkov-1974}, 1974\\
 $4$ & ${\mathsf{H}_4}$    &  $1.3625\ldots$  &  Horv\'ath \cite{horvath-1980}, 1980\\
 $5$ & ${\mathsf{H}_5}$    &  $1.4494\ldots$  &  Horv\'ath \cite{horvath-1986}, 1986\\
\hline
\end{tabular} 
\medskip
\caption{Optimal lattice packing-coverings.}
\end{center}
\end{table}

\section{Lattices and Positive Quadratic Forms}
\label{sec:lat-quad}

It is sometimes convenient to switch from the language of lattices to
the language of positive definite quadratic forms (PQFs from now
on). In this section we give a dictionary. For further reading we
refer to \cite{cs-1988} and \cite{sv-2004}.

Given a $d$-dimensional lattice $L=A\Z^d$ with basis $A$ we associate
a $d$-dimensional PQF $ Q[\vec{x}] =\vec{x}^t A^t A \vec{x}
=\vec{x}^t G \vec{x} $, where the \textit{Gram matrix} $G = A^t A$ is
symmetric and positive definite.  We will carelessly identify
quadratic forms with symmetric matrices by saying $Q = G$ and
$Q[\vec{x}] = \vec{x}^t Q \vec{x}$.  The set of quadratic forms is a
$\binom{d+1}{2}$-dimensional real vector space $\sd$, in which the
set of PQFs forms an open, convex cone~$\sdo$.  The PQF~$Q$ depends on
the chosen basis $A$ of $L$. For two arbitrary bases $A$ and $B$ of
$L$ there exists a $U\in\gldz$ with $A = BU$.  Thus, $\gldz$ acts on
$\sdo$ by $ Q \mapsto U^t Q U$. A PQF $Q$ can be associated to different
lattices $L = A\Z^d$ and $L' = A'\Z^d$. In this case there exists an
orthogonal transformation $O$ with $A = OA'$. Note that the packing
and covering density, as well as the packing-covering constant, 
are invariant with respect to orthogonal transformations.

The \textit{determinant} (or discriminant) of a PQF $Q$ is defined by
$\det(Q)$.  The \textit{homogeneous minimum} $\lambda(Q)$ and the
\textit{inhomogeneous minimum} $\mu(Q)$ are given by
\[
\lambda(Q) = \min_{\vec{v} \in \Z^d \setminus \{\vec{0}\}} Q[\vec{v}],
\qquad
\mu(Q) = \max_{\vec{x} \in \R^d} \min_{\vec{v} \in \Z^d} Q[\vec{x} - \vec{v}].
\]
If $Q$ is associated to $L$, then $\det(L) = \sqrt{\det(Q)}$,
$\mu(L) = \sqrt{\mu(Q)}$, $\lambda(L) = \sqrt{\lambda(Q)}/2$.
Using this dictionary we define
\begin{eqnarray*}
\delta(Q) = \sqrt{\frac{\lambda(Q)^d}{\det Q}} \cdot \frac{\kappa_d}{2^d},
\;\;
\Theta(Q) = \sqrt{\frac{\mu(Q)^d}{\det Q}} \cdot \kappa_d
\;\;\mbox{and}\;\;
\gamma(Q) = 2\sqrt{\frac{\mu(Q)}{\lambda(Q)}}
.
\end{eqnarray*}

We say that a lattice $L$ with associated PQF $Q$ gives a
\textit{locally optimal lattice packing},
\textit{locally optimal lattice covering} or
\textit{locally optimal lattice packing-covering},
if there is a neighborhood of $Q$ in $\sdo$, so
that we have 
$\delta(Q) \geq \delta(Q')$, 
$\Theta(Q) \leq \Theta(Q')$ or 
$\gamma(Q) \leq \gamma(Q')$
respectively, for all $Q'$ in this neighborhood.

\section{On Packing Lattices}
\label{sec:pack-lat}

A PQF $Q$ attaining a local maximum of $\delta(Q)$ is
called \textit{extreme}. A PQF attaining $\delta_d = \max_{Q \in \sdo}
\delta(Q)$ is called \textit{absolutely extreme} or \textit{critical}.
One can characterize an extreme PQF using the geometry of its minimal
vectors $$\Min(Q) = \{\vec{v} \in \Z^d : Q[\vec{v}] = \lambda(Q)\}.$$  
Before we state the characterization in Theorem~\ref{th:extremevoro}, we give some more definitions.
A PQF $Q'$ is called \textit{perfect} if it is uniquely determined by
its minimal vectors, i.e.\ $Q'$ is the unique solution of the linear
equations $Q[\vec{v}] = \lambda(Q')$, $\vec{v} \in \Min(Q')$.  A PQF
$Q$ is called \textit{eutactic} if 
\[
Q^{-1}
\in\relint\cone\{ \vec{v}\vec{v}^t : \vec{v}\in \Min(Q)\}
= \left\{\sum_{\vec{v}\in \Min(Q)} \lambda_{\vec{v}} \vec{v}\vec{v}^t: \lambda_{\vec{v}} > 0, \vec{v}\in \Min(Q)\right\}
.
\]
The eutaxy and the polyhedral cone $\cone\{ \vec{v}\vec{v}^t : \vec{v}\in \Min(Q)\}$ also plays an important role 
in the lattice packing-covering problem.
As a general reference on basic facts about polyhedral cones, 
which are used throughout this article, 
we refer to the book of Ziegler \cite{ziegler-1998}.

\begin{theorem}[Voronoi \cite{voronoi-1907}]
\label{th:extremevoro} 
A PQF is extreme if and only if it is perfect and eutactic.
\end{theorem}

This provides an easy way for proving that a given PQF is extreme:
after finding the minimal vectors, one has to solve a system of linear
equations to show its perfectness. Then, one has to solve a linear
programming problem to verify its eutaxy. By scaling we can
normalize an extreme PQF $Q$ so that $\lambda(Q)$ is rational. Then,
since $Q$ is perfect, the matrix entries of $Q$ are rational as well.

It turns out that there exist only finitely many pairwise
non-equivalent perfect PQFs in $\sdo$. We want to describe two
classical algorithms to attain all perfect forms of a given dimension
$d$. The first one goes back to Minkowski, the second one is due to
Voronoi. Here we only state definitions and main
results. Additionally, we briefly sketch the computations which we
were able to perform in low dimensions. We compare them with
corresponding results in the literature.

For history and further remarks we refer to \cite{gl-1987},
\cite{rb-1979}, \cite{vdwaerden-1956}, \cite{martinet-2003} and to
references therein. Chapter~\S v of Gruber and Lekkerkerker's
book~\cite{gl-1987} gives a comprehensive survey about history,
results and literature of the reduction theory of PQFs. The
article~\cite{rb-1979} introduces to methods for studying the geometry
of PQFs and contains many proofs. Van der Waerden's
paper~\cite{vdwaerden-1956} is a classic resource for Minkowski's
approach. The recent book~\cite{martinet-2003} of Martinet gives a
contemporary view on Voronoi's approach and on possible
generalizations.

\subsection{Minkowski's Approach}

\begin{definition}
A PQF $Q =(q_{ij}) \in \sdo$ is called \textit{Minkowski reduced} if
\begin{enumerate}
\item[(i)] $Q[\vec{x}] \geq q_{ii}$ whenever $\gcd(\vec{x}_i, \ldots,
\vec{x}_d)=1$, $i = 1, \ldots, d$,
\item[(ii)] $q_{i,i+1} \geq 0$, $i = 1, \ldots, d-1$.
\end{enumerate}
\end{definition}

Every PQF is equivalent to a Minkowski reduced PQF. The following
procedure, which is nothing but an algorithmic interpretation of the
definition, finds a Minkowski reduced PQF equivalent to a given PQF
$Q$. Choose a minimal vector $\vec{v}_1$ of $Q$. Then, choose among
all vectors in $\Z^d$, which can complement $\vec{v}_1$ to a lattice
basis of $\Z^d$, a vector $\vec{v}_2$ for which the value
$Q[\vec{v}_2]$ is minimal. Using this greedy strategy we get a basis
$V = (\vec{v}_1, \ldots, \vec{v}_d) \in \gldz$. Now we choose signs
for $\vec{v}_i$ so that $\vec{v}_{i-1}^t Q \vec{v}_i \geq 0$. Hence,
$V^t Q V$ is Minkowski reduced.

The set of Minkowski reduced PQFs forms an unbounded cone in $\sdgeo$
which is defined by the linear inequalities (i) and (ii). By $\MM$ 
we denote the cone defined by the linear inequalities (i) and by
$\MM^+$ the one which is defined by the linear inequalities (i)
and (ii). Minkowski \cite{minkowski-1905} showed that $\MM$, and
hence $\MM^+$, is a polyhedral cone, i.e.\ that finitely many
inequalities (i) imply all others. He showed that every extreme PQF is
equivalent to a PQF lying on a ray (a one-dimensional face) of $\MM^+$. 
Ryshkov \cite{ryshkov-1970} proved that every perfect PQF is
equivalent to a PQF lying on a ray of $\MM^+$. On the other hand
Cohn, Lomakina and Ryshkov \cite{clr-1980} found a ray of $\MM^+
\subset {\MS}^5_{\geq 0}$ which contains non-perfect PQFs.

Minkowski \cite{minkowski-1887} gave a list of conditions implying all
others in (i) up to dimension~$6$. Tammela \cite{tammela-1977}
enlarged this list to dimension $7$. Besides the $d - 1$ inequalities
$q_{11} \leq \ldots \leq q_{dd}$, the linear conditions for $\MM$, 
$d = 2, \ldots, 7$, are attained by plugging the values
$\vec{x} = (x_1, \ldots, x_d) \in \Z^d$ from Table~\ref{tammela-tab}
into (i), where the indices $i_1, \dots, i_d$ 
run through all permutations of $\{1, \dots, d\}$.  
If $d < 7$ one has to omit the columns $d + 1,
\dots, 7$ and the rows with more than $d$ non-zero entries.

\begin{table}[ht]
\begin{center}
\begin{tabular}{cccccccc}
$\;\;x_{i_1}$ & $\pm x_{i_2}$ & $\pm x_{i_3}$ & $\pm x_{i_4}$ & 
            $\pm x_{i_5}$ & $\pm x_{i_6}$ & $\pm x_{i_7}$ & \\
\hline\hline
1 & 1 & 0 & 0 & 0 & 0 & 0 & \\
1 & 1 & 1 & 0 & 0 & 0 & 0 & \\
1 & 1 & 1 & 1 & 0 & 0 & 0 & \\
1 & 1 & 1 & 1 & 1 & 0 & 0 & \\
2 & 1 & 1 & 1 & 1 & 0 & 0 & \\
1 & 1 & 1 & 1 & 1 & 1 & 0 & \\
2 & 1 & 1 & 1 & 1 & 1 & 0 & \\
2 & 2 & 1 & 1 & 1 & 1 & 0 & \\
3 & 2 & 1 & 1 & 1 & 1 & 0 & $^1$ \\
1 & 1 & 1 & 1 & 1 & 1 & 1 & \\
2 & 1 & 1 & 1 & 1 & 1 & 1 & \\ 
3 & 1 & 1 & 1 & 1 & 1 & 1 & \\ 
2 & 2 & 1 & 1 & 1 & 1 & 1 & \\ 
3 & 2 & 1 & 1 & 1 & 1 & 1 & \\ 
2 & 2 & 2 & 1 & 1 & 1 & 1 & \\ 
3 & 2 & 2 & 1 & 1 & 1 & 1 & \\ 
4 & 2 & 2 & 1 & 1 & 1 & 1 & \\ 
3 & 3 & 2 & 1 & 1 & 1 & 1 & $^2$ \\   
4 & 3 & 2 & 1 & 1 & 1 & 1 & $^3$ \\ 
3 & 2 & 2 & 2 & 1 & 1 & 1 & \\ 
4 & 3 & 2 & 2 & 1 & 1 & 1 & \\
\hline
\end{tabular}
\medskip
\caption{Tammela's list of linear inequalities defining $\MM$}
\label{tammela-tab}
\end{center}
\end{table}

We checked these conditions for redundancy using the software
\texttt{lrs} of Avis \cite{avis-2004}. Our computations show that in
row $^1$, for $d=6$, the entries with $x_1 \neq 3$ 
are redundant (as already mentioned by Tammela in \cite{tammela-1973}).
For $d=7$ redundant entries are those with $x_1 \neq 0,3$
and $x_1=0, x_2 \neq 3$.
In row $^2$ the entries with $x_1 \neq 3$ or
$x_2 \neq 3$, and in row $^3$ the entries with $x_1 \neq 4$ or $x_2
\neq 3$ are redundant.  The remaining conditions are all non-redundant
and define a facet of $\MM$.  In Table~\ref{fac-rays-numbers-tab} we
list the number of facets and rays as far as we were able to compute
them with \texttt{cdd} \cite{fukuda-2003}. We hereby confirm earlier
results by ${}^a$ Minkowski \cite{minkowski-1887} and ${}^b$ Barnes
and Cohn \cite{bk-1976}. We made the data available from the
\href{www.arxiv.org}{arXiv.org} e-print archive. To access it,
download the source files for the paper
\href{http://www.arxiv.org/math/0412320}{arXiv:math.MG/0412320}.  The
files \texttt{mink2.ine}, \dots, \texttt{mink6.ine} (due to its size the file
\texttt{mink7.ine} is only available from the authors) contain the facets of
$\MM$, the files \texttt{mink2.ext}, \dots, \texttt{mink6.ext} contain
the rays of $\MM$, the files \texttt{minkp2.ine}, \dots,
\texttt{minkp6.ine} contain the facets of $\MM^+$, and the files
\texttt{minkp2.ext}, \dots, \texttt{minkp5.ext} contain the rays of
$\MM^+$.  For the data format we chose the common convention
(Polyhedra $H$-format for the \texttt{*.ine}-files and Polyhedra
$V$-format for the \texttt{*.ext}-files) of the software packages
\texttt{cdd} and \texttt{lrs}.

The computational bottlenecks of Minkowski's approach are apparent. It
is not easy to find a sufficiently small system of linear inequalities
defining $\MM$ (or of $\MM^+$). Even if one has a minimal system of
linear inequalities, then computing its rays is a very difficult
computational problem in higher dimensions.

\begin{table}[ht]
\begin{center}
\begin{tabular}{|c|c|c|c|c|c|}
\hline 
$d$ & $\binom{d+1}{2}$ & \# Facets $\MM$ & \# Rays $\MM$  & \#
Facets $\MM^+$  & \# Rays $\MM^+$  \\
\hline\hline 
$2$ & $3$ & $3^a$ & $3^a$ & $3^a$ & $3^a$ \\
$3$ & $6$ & $12^a$ & $19^b$ & $9^a$ & $11^b$ \\
$4$ & $10$ & $39^a$ & $323^b$ & $26^b$ & $109^b$ \\
$5$ & $15$ & $200$ & $15971$ & $117$ & $4105$ \\
$6$ & $21$ & $1675$ & ? & $1086$ &  ?\\
$7$ & $28$ & $65684$ & ? & ? & ? \\
\hline
\end{tabular}
\medskip
\caption{Known number of facets and rays of $\MM$ and $\MM^+$.}
\label{fac-rays-numbers-tab}
\end{center}
\end{table}

\subsection{Voronoi's Approach}

Now we describe Voronoi's algorithm \cite{voronoi-1907} for finding all
perfect forms of a given dimension. Let $m$ be a positive number. In
the remaining of this section we assume that every perfect form $Q$ is scaled so that
$\lambda(Q) = m$. The set
\[
\MP_m = 
\left\{ Q \in \sdo : \mbox{$Q[\vec{v}] \geq m$ for all $\vec{v} \in \Z^d \setminus \{0\}$} \right\}
\]
is a convex, locally finite polyhedral cone. Its boundary consists of
the PQFs with homogeneous minimum~$m$. A PQF $Q$ is perfect if and
only if it is a vertex of $\MP_{m}$. The set of perfect PQFs of a
given dimension $d$ naturally carries a graph structure which we
denote as the \textit{Voronoi graph in dimension $d$}: Two perfect
PQFs $Q$, $Q'$ are connected by an edge if the line segment
$\conv\{Q,Q'\}$ (convex hull of $Q$ and $Q'$)
is an edge of~$\MP_m$. In this case we say that $Q$
and $Q'$ are \textit{Voronoi neighbors}. The group $\mbox{GL}_d(\Z)$
acts on $\MP_m$, on its vertices and on its edges by $Q \mapsto U^t Q
U$. Therefore, one can enumerate perfect PQFs by a graph traversal
algorithm which we shall now describe.

Voronoi's first perfect form $Q[\vec{x}] = \sum_{i = 1}^d x_i^2 + \sum_{i < j}
x_i x_j$, which is associated to the root lattice $\mathsf{A}_d$, can
serve as a starting point in any dimension. For implementing a graph
traversal algorithm one has to find the Voronoi neighbors of a given
perfect form~$Q$. Consider the unbounded polyhedral cone
\[
\MP(Q) = \{Q'\in\sd : Q'[\vec{v}] \geq m, \vec{v} \in \Min(Q)\}.
\]
We compute the rays $Q + \R_{\geq 0}R_i$, $i = 1,\ldots,n$ of 
$\MP(Q)$. The $R_i$ turn out to be indefinite quadratic forms. So
there are $\vec{v} \in \Z^d$ with $R_i[\vec{v}] < 0$. Then, the
Voronoi neighbors of $Q$ are $Q + \rho_i R_i$ where $\rho_i$ is the
smallest positive number so that $\lambda(Q + \rho_i R_i) = m$ and
$\Min(Q + \rho_i R_i) \not \subseteq \Min(Q)$. 
It is possible to determine $\rho_i$, for example with 
the following procedure:

\begin{center}
\begin{minipage}{8cm}
\begin{flushleft}
$(\alpha, \beta) \leftarrow (0,1)$\\
\smallskip
\textbf{while} $Q + \beta R_i \not\in \sdo$ or $\lambda(Q + \beta R_i) = m$ \textbf{do}\\
\hspace{2ex} \textbf{if} $Q + \beta R_i \not\in \sdo$ \textbf{then} $\beta \leftarrow \frac{\alpha + \beta}{2}$\\
\hspace{2ex} \textbf{else} $(\alpha,\beta) \leftarrow (\beta, 2\beta)$\\
\hspace{2ex} \textbf{end if}\\
\textbf{end while}\\
\smallskip
\textbf{while} $\Min(Q + \alpha R_i) \subseteq \Min Q$ \textbf{do}\\
\hspace{2ex} $\gamma \leftarrow \frac{\alpha + \beta}{2}$\\
\hspace{2ex} \textbf{if} $\lambda(Q + \gamma R_i) < m$ \textbf{then} $\beta \leftarrow \gamma$\\
\hspace{2ex} \textbf{else} $\alpha \leftarrow \gamma$\\
\hspace{2ex} \textbf{end if}\\
\textbf{end while}\\
\smallskip
$\rho_i \leftarrow \alpha$
\end{flushleft}
\end{minipage}
\end{center}

Perfect forms were classified up to dimension $7$. A list of all these
forms is given in the paper \cite{cs-1988b} of Conway and Sloane. One
can find an electronic version in the Catalogue of Lattices\footnote{
\href{http://www.research.att.com/~njas/lattices/perfect.html}{http://www.research.att.com/\~{}njas/lattices/perfect.html}
} 
by Nebe and Sloane.
On his homepage\footnote{
\href{http://www.math.u-bordeaux.fr/~martinet}{http://www.math.u-bordeaux.fr/\~{}martinet}
}
,
Martinet reports that up to now, $10916$ pairwise inequivalent 
perfect forms are known in dimension $8$ and lists them.

We verified the results for dimensions $d\leq 6$
using the programs \texttt{lrs} by Avis
\cite{avis-2004}, \texttt{isom} by Plesken and Souvignier
\cite{ps-1997} and \texttt{shvec} by Vallentin
\cite{vallentin-1999}. In Table~\ref{perfect-tab} we give the known
classifications of perfect forms, extreme forms and absolute extreme
forms together with the references where the classifications were
established.

\begin{center}
\begin{table}[ht]
\begin{tabular}{|c|c|c|c|c|}
\hline 
$\; \mathbf{d} \;$ & $\; \mathbf{\binom{d+1}{2}} \;$ & \textbf{\# perfect} & 
                    \textbf{\# extreme} & \textbf{\# critical} \\
\hline 
$2$ & $3$ & $1^a$ & $1^a$ & $1^a$ \\
$3$ & $6$ & $1^b$ & $1^b$ & $1^b$ \\
$4$ & $10$ & $2^c$ & $2^c$ & $1^c$ \\
$5$ & $15$ & $3^d$ & $3^d$ & $1^d$ \\
$6$ & $21$ & $7^f$ & $6^e$ & $1^f$ \\
$7$ & $28$ & $33^h$ & $30^h$ & $1^g$ \\
$8$ & $36$ & $\geq 10916^i$ & ? & $1^g$ \\
$9$ & $45$ & $> 500,000$ & ? &  ?\\[0.3cm]
$24$ & $300$ & ? &  ? & $1^j$\\
\hline
\end{tabular}

\begin{tabular}{ll}
${}^a$ & Lagrange \cite{lagrange-1773}\\
${}^b$ & Gau{\ss} \cite{gauss-1840}\\
${}^c$ & Korkine, Zolotareff \cite{kz-1873}\\
${}^d$ & Korkine, Zolotareff \cite{kz-1877}\\
${}^e$ & Hofreiter \cite{hofreiter-1933}\\
${}^f$ & Barnes \cite{barnes-1957a}\\
${}^g$ & Vetchinkin \cite{vetchinkin-1980}\\
${}^h$ & Jaquet-Chiffelle \cite{jaquet-1993}\\
${}^i$ & Laihem, Baril, Napias, Batut, Martinet (see \href{http://www.math.u-bordeaux.fr/~martinet}{http://www.math.u-bordeaux.fr/\~{}martinet})\\
${}^j$ & Cohn, Kumar \cite{ck-2004}
\end{tabular}
\medskip
\caption{Classifications of perfect, extreme and absolute
extreme forms}
\label{perfect-tab}
\end{table}
\end{center}

The computational bottleneck of Voronoi's approach is mainly the
enumeration of all rays of the polyhedral cone $\MP(Q)$ in case 
of a large set $\Min(Q)$ of minimal vectors. Martinet writes in
\cite{martinet-2003}, Ch.\ 7.11 : "The existence of $\mathsf{E}_8$ [...]
makes hopeless any attempt to construct the Voronoi graph in dimension
$8$". Another problem is the combinatorial explosion, when $d\geq 9$.
We found more than $500000$ inequivalent perfect forms in dimension $9$
and we strongly believe there exist millions of them.

Finally, we want to remind of Coxeter's $\mathsf{A}_d$-hypothesis: Although
finding perfect forms with maximal packing density is a very difficult
problem, finding perfect forms with minimal packing density might be
very easy. In \cite{coxeter-1951} Coxeter formulates the following
conjecture:

\begin{conjecture} (Coxeter's $\mathsf{A}_d$-hypothesis)
Voronoi's first perfect form gives the minimal packing density among
all perfect forms of a given dimension.
\end{conjecture} 

Our computations support Coxeter's conjecture. But on the contrary,
Conway and Sloane \cite{cs-1988b} conjecture that it is false for
sufficiently large $d$.

\section{Covering and Packing-Covering Lattices}
\label{sec:cov-lat}

The lattice covering problem and the lattice packing-covering problem
can be treated in parallel.
We describe below that both problems have only finitely many local optima
which can be found by solving finitely many convex optimization problems.
This is mainly due to Voronoi's theory of Delone subdivisons, 
which we briefly review. For a detailed account we refer to \cite{sv-2004}.

With an implementation of the proposed algorithms
we found all local optima in dimension $d\leq 5$ and
some new best known lattices in dimension $d\geq 6$.

For both problems, 
recognition of local optima is not as easy as for the lattice packing problem.
Due to the involved convexity we can give sufficient conditions
for local optima, allowing to compute a certificate for the local optimality
of a lattice. This is in particular applicable, if the Delone subdivison is
a triangulation which is the generic case (for definitions see below). 
Exemplarily we give a proof of the local packing-covering 
optimality of the lattices $\mathsf{A}^*_d$.

In some cases it is possible to attain good or even tight ``local lower bounds'' for
the lattice covering density and the packing-covering constant.
This is demonstrated for the local packing-covering optimality of the Leech lattice.
A similar proof of the local covering optimality of the Leech lattice is given in
\cite{sv-2004b}.

\subsection{Voronoi's Theory of Delone Subdivisions}

Let $Q$ be a positive semidefinite quadratic form.  A polyhedron $P =
\conv\{\vec{v}_1, \vec{v}_2, \ldots\}$ with $\vec{v}_1, \vec{v}_2,
\ldots \in \Z^d$, is called a \textit{Delone polyhedron} of~$Q$ if
there exists a $\vec{c} \in \R^d$ and a real number $r \in \R$ with
$Q[\vec{v}_i - \vec{c}] = r^2$ for all $i = 1, 2, \ldots$, and
$Q[\vec{v} - \vec{c}] > r^2$ for all other $\vec{v} \in \Z^d \setminus
\{\vec{v}_1, \vec{v}_2, \ldots\}$.  The set $\Del(Q)$ of all Delone
polyhedra is called the \textit{Delone subdivision} of~$Q$. It is a
periodic face-to-face tiling of $\R^d$. Therefore $\Del(Q)$ is
completely determined by all Delone polytopes having a vertex at the
origin $\vec{0}$. We call two Delone polyhedra $L,L'$
\textit{equivalent} if there exists a $\vec{v} \in \Z^d$ so that $L =
\vec{v} \pm L'$.  Note moreover that the inhomogeneous minimum
$\mu(Q)$ is at the same time the maximum squared circumradius of its
Delone polyhedra. We say that the Delone subdivision of a positive
semidefinite quadratic form $Q'$ is a \textit{refinement} of the
Delone subdivision of $Q$, if every Delone polytope of $Q'$ is
contained in a Delone polytope of $Q$.

By a theory of Voronoi \cite{voronoi-1908}, the
set of positive semidefinite quadratic forms with a fixed Delone
subdivision $\MD$ is an open (with respect to its affine hull)
polyhedral cone in $\sdgeo$.  We refer to this set as the
\textit{secondary cone} $\SC(\MD)$ of the subdivision.  In the
literature the secondary cone is sometimes called $L$-type domain of
the subdivision. The topological closure $\overline{\SC(\MD)}$
of a secondary cone is a closed polyhedral cone. 
The relative interior of each face in $\sdo$ is the secondary cone 
of another Delone subdivision. If a face is contained in the boundary 
of a second face, then the corresponding Delone subdivision of the 
first is a true refinement of the second one.

The interior of faces of maximal dimension $\binom{d+1}{2}$ contain
PQFs whose Delone subdivision is a triangulation, that is, it consists
of simplices only. We refer to such a subdivision as a
\textit{simplicial} Delone subdivision or \textit{Delone
triangulation}.  As mentioned in Section~\ref{sec:lat-quad}, the group
$\GL_d(\Z)$ acts on $\sdgeo$.  One of the key observations of Voronoi
is that under this group action there exist only finitely many
inequivalent Delone subdivisions, respectively secondary cones.

\begin{theorem}[Voronoi \cite{voronoi-1908}]
\label{th:mainvoronoi} 
The topological closures of secondary cones of Delone
triangulations give a face-to-face tiling of $\sdgeo$.    
The group $\GL_d(\Z)$ acts on the tiling, 
and under this group action there are only finitely  
many non-equivalent secondary cones.
\end{theorem}

Given a Delone triangulation $\MD$, the Delone triangulations $\MD'$
with $\overline{\SC(\MD')}$ sharing a facet with $\overline{\SC(\MD)}$
are attained by {\em bistellar operations} ({\em flips}). These change a
triangulation only in certain {\em repartitioning polytopes} associated to the
facet.  By this operation it becomes possible to enumerate all Delone
triangulations, and hence all Delone subdivisions in a given
dimension. For details we refer to \cite{sv-2004}.

\subsection{Obtaining Local Optima via Convex Optimization}

For a fixed triangulation $\MD$, we can formulate the lattice
covering, as well as the lattice packing-covering problem in the
framework of {\em Determinant Maximization Problems}.  Following
Vandenberghe, Boyd, and Wu \cite{vbw-1998} their general form is
\begin{equation}   \label{eq:maxdet-problem}
\begin{array}{ll}
\mbox{\textbf{minimize}} & \vec{c}^t\vec{x}-\log\det G(\vec{x})\\
\mbox{\textbf{subject to}} & G(\vec{x}) \succ 0, 
                             F(\vec{x}) \succeq 0.
\end{array}
\end{equation} 
Here, the optimization vector is $\vec{x} \in \R^D$.  The objective
function contains a linear part given by $\vec{c} \in \R^D$ and $G :
\R^D \to \R^{m \times m}$, $F : \R^D \to \R^{n \times n}$ are both
affine maps
$$
\begin{array}{rcl}
G(\vec{x}) & = & G_0 + x_1 G_1 + \cdots + x_D G_D,\\
F(\vec{x}) & = & F_0 + x_1 F_1 + \cdots + x_D F_D,
\end{array}
$$ 
where $G_i \in \R^{m \times m}$, $F_i \in \R^{n \times n}$, $i =
0,\ldots, D$, are symmetric matrices.  The notation $G(\vec{x}) \succ
0$ and $F(\vec{x}) \succeq 0$ gives the constraints ``$G(\vec{x})$ is
positive definite'' and ``$F(\vec{x})$ is positive semidefinite''.
Note that we are dealing with a
so-called \textit{semidefinite programming problem}, if $G(\vec{x})$
is the identity matrix for all $\vec{x} \in \R^D$.

For the lattice covering, as well as the lattice packing-covering
problem, we can express $\mu(Q)\leq 1$ 
as a {\em linear matrix inequality} $F(q_{ij}) \succeq 0$
with optimization vector $Q = (q_{ij})$.
To see this, it is crucial to observe that
an inner product $(\cdot, \cdot)$ defined by 
$(\vec{x}, \vec{y}) = \vec{x}^tQ\vec{y}$
gives a linear expression in the parameters $(q_{ij})$
for any fixed choice $\vec{x},\vec{y}\in \Z^d$ (or $\R^d$).
Delone, Dolbilin, Ryshkov and Stogrin \cite{ddrs-1970}
showed 

\begin{proposition}
\label{prop:lmi}
Let $L = \conv\{\vec{0}, \vec{v}_1, \ldots, \vec{v}_d\}\subseteq
\R^d$ be a $d$-dimensional simplex. 
Then $L$'s circumradius is at most~$1$ with respect to $(\cdot, \cdot)$ 
if and only if
$$
\BR_L(Q) =
\begin{pmatrix}
4 & (\vec{v}_1, \vec{v}_1) & (\vec{v}_2, \vec{v}_2) & \ldots &
(\vec{v}_d, \vec{v}_d)\\
(\vec{v}_1, \vec{v}_1) & (\vec{v}_1, \vec{v}_1) & (\vec{v}_1, \vec{v}_2) & \ldots &
(\vec{v}_1, \vec{v}_d)\\
(\vec{v}_2, \vec{v}_2) & (\vec{v}_2, \vec{v}_1) & (\vec{v}_2,
\vec{v}_2) & \ldots & (\vec{v}_2, \vec{v}_d)\\
\vdots & \vdots & \vdots & \ddots & \vdots\\
(\vec{v}_d, \vec{v}_d) & (\vec{v}_d, \vec{v}_1) & (\vec{v}_d, \vec{v}_2) & \ldots &
(\vec{v}_d, \vec{v}_d)\\
\end{pmatrix}
\succeq 0.
$$
\end{proposition}

Since a block matrix is semidefinite if and only if the blocks are
semidefinite, we conclude

\begin{proposition}  
\label{propos:mu}
Let $Q = (q_{ij}) \in \sdo$ be a PQF. Let $\MD$ be a Delone
triangulation refining $\Del(Q)$, and let $L_1, \ldots, L_n$ be a
representative system of all non-equivalent $d$-dimensional Delone
polytopes in $\MD$. Then
$$
\mu(Q)\leq 1
\quad\Longleftrightarrow\quad
\begin{pmatrix}
\boxed{\BR_{L_1}(Q)} & 0                   & 0 & \ldots & 0\\
0                   & \boxed{\BR_{L_2}(Q)} & 0 & \ldots & 0\\
\hdotsfor{5}\\
0                   & 0                   & 0 & \ldots & \boxed{\BR_{L_n}(Q)}
\end{pmatrix}
\succeq 0.
$$
\end{proposition}

Thus $\mu(Q)\leq 1$ can be brought into one linear matrix inequality
of type $F(q_{ij}) \succeq 0$.  We can moreover add linear constraints
on the parameters $q_{ij}$ by extending $F$ by a $1\times 1$ block
matrix for each linear inequality. In this way we can get one linear
matrix inequality for the two constraints $\mu(Q)\leq 1$ and $Q \in
\overline{\SC(\MD)}$.

For a fixed Delone triangulation~$\MD$, we can therefore determine the
optimal solutions of the lattice covering problem of all PQFs $Q$ for
which $\MD$ is a refinement of~$\Del(Q)$.  Recall that the covering
density of a PQF~$Q$ in $d$~variables is $\Theta(Q) =
\sqrt{\frac{\mu(Q)^d}{\det Q}} \cdot \kappa_d$.  Scaling of $Q$ by a
positive real number~$\alpha$ leaves $\Theta$ invariant. Thus we may
maximize $\det(Q)$ while $\mu(Q) \leq 1$. For all $Q\in
\overline{\SC(\MD)}$ this can be achieved by solving
$$
\begin{array}{ll}
\mbox{\textbf{minimize}} & -\log\det(Q)\\
\mbox{\textbf{subject to}} & \mbox{$Q\succ 0$,}\\ 
         & \mbox{$Q \in \overline{\SC(\MD)}$, $\mu(Q) \leq 1$}.\\
\end{array}
$$ 

With an analogues specialization of problem \eqref{eq:maxdet-problem}, we
are able to attain optimal solutions of the lattice packing-covering
problem among all PQFs $Q$ for which $\MD$ is a refinement of~$\Del(Q)$.
Because $\gamma(Q) = 2 \cdot\sqrt{\mu(Q)/\lambda(Q)}$, we have to
maximize $\lambda(Q)$ while $\mu(Q)\leq 1$.  Maximizing $\lambda(Q)$
is not as straightforward as maximizing the determinant, since we do
not know which vector is the shortest.  By a theorem of Voronoi
\cite{voronoi-1908} we know though that among the (at most $2(2^d -
1)$) edges $[\vec{0}, \vec{v}] = \conv\{\vec{0}, \vec{v}\}\in \MD$,
there exist some with $\vec{v}\in\Min(Q)$.  Consequently, if we
require $Q[\vec{v}] \leq Q[\vec{w}]$ for all $\vec{w}$ with $[\vec{0},
\vec{w}]\in \MD$, we know $\vec{v}\in\Min(Q)$, respectively
$\lambda(Q)=Q[\vec{v}]$.  So we can solve for each 
$[\vec{0}, \vec{v}]\in \MD$ the semidefinite programming problem

$$
\begin{array}{ll}
\mbox{\textbf{minimize}} & -Q[\vec{v}]\\
\mbox{\textbf{subject to}} & \mbox{$Q[\vec{v}] \leq Q[\vec{w}]$, for all $\vec{w}$ with $[\vec{0}, \vec{w}]\in \MD$,}\\
& \mbox{$Q \in \overline{\SC(\MD)}$}, \mbox{$\mu(Q) \leq 1$}.\\
\end{array}
$$

Note that in many cases the constraints have no feasible solutions,
since in general not all of the $\vec{v}$ with $[\vec{0}, \vec{v}]\in \MD$ 
are elements of $\Min(Q)$.

So both, the lattice covering as well as the lattice packing-covering problem,
are reduced to convex programming problems if restricted to the closure
of a secondary cone.
Consequently, there is at most one local minimum of 
$\Theta$, respectively $\gamma$, for each of the cones $\overline{\SC(\MD)}$.
This was first observed by Barnes and Dickson
\cite{bd-1967} in the covering case and by Ryshkov \cite{ryshkov-1974}
for the packing-covering problem.

\begin{proposition} 
\label{prop:invariance}
Let $\MD$ be a Delone triangulation. Then there exists a unique minimum of 
\begin{enumerate}
\item \label{item:local-opt:1}
$\{ \Theta(Q) : Q \in \overline{\SC(\MD)}\}$ and the set of PQFs attaining it 
is equal to all positive multiples of a single PQF.  
\item 
$\{ \gamma(Q) : Q \in \overline{\SC(\MD)}\}$ and the set of PQFs attaining it 
is convex.
\end{enumerate}
\end{proposition}

In case of a triangulation $\MD$,
it follows that $Q\in \SC(\MD)$ is a locally optimal solution 
with respect to $\Theta$ or $\gamma$ if and only if
it is an optimal solution within $\overline{\SC(\MD)}$.
In general we have the trivial

\begin{proposition}
A PQF $Q$ is a locally optimal solution with respect to $\Theta$ or
$\gamma$, if and only if it is an optimal solution for all
Delone triangulations $\MD$ with $Q\in\overline{\SC(\MD)}$.
\end{proposition}

\subsection{Computational Results}

By the foregoing propositions we know that the number of local optima
is bounded from above by the number of pairwise inequivalent Delone
triangulations in $\R^d$.  Voronoi \cite{voronoi-1908} classified
these triangulation in dimension $2$, $3$ (only one each) and $4$
(three).  By the work of Baranovskii and Ryshkov \cite{br-1973},
\cite{rb-1976} Engel \cite{engel-1998}, and Engel and Grishukhin
\cite{eg-2002} we know of exactly $222$ Delone triangulations in
dimension $d=5$.  Using \texttt{lrs} \cite{avis-2004} and an
implementation (in C++) of Voronoi's algorithm for enumerating
Delone triangulations, we were able to confirm these results
\cite{sv-2004}.  For dimension $d\geq 6$ we experience a combinatorial
explosion, e.g.\ Engel \cite{engel-2004} reports on more than
$2,129,120$ pairwise inequivalent Delone triangulations for $d=6$.

\begin{center}
\begin{table}[ht]
\begin{tabular}{|c|c|c|c|}
\hline 
$\; \mathbf{d} \;$ & $\; \mathbf{\binom{d+1}{2}} \;$ & \textbf{\# covering optima} & 
                    \textbf{\# packing-covering optima}  \\
\hline 
$2$ & $3$ & $1^a$ & $1^b$  \\
$3$ & $6$ & $1^c$ & $1^b$ \\
$4$ & $10$ & $3^d$ & $3^e$ \\
$5$ & $15$ & $222^f$ & $47^f$ \\
\hline
\end{tabular}

\begin{tabular}{ll}
${}^a$ & Kershner \cite{kershner-1939}\\
${}^b$ & Ryshkov \cite{ryshkov-1974}\\
${}^c$ & Bambah \cite{bambah-1954a}\\
${}^d$ & Baranovski \cite{baranovskii-1965}, Dickson \cite{dickson-1967}\\
${}^e$ & Horv\'ath \cite{horvath-1980} \\
${}^f$ & Sch\"urmann, Vallentin \cite{sv-2004}
\end{tabular}
\medskip
\caption{Classifications of locally optimal covering and packing-covering lattices.}
\label{table:latpackcov}
\end{table}
\end{center}

For each of the triangulations in dimension $d\leq 5$ we
determined the local optima with respect to $\Theta$ and $\gamma$
using the software package
\texttt{MAXDET}\footnote{\href{http://www.stanford.edu/~boyd/MAXDET.html}{http://www.stanford.edu/\~{}boyd/MAXDET.html}}
of Wu, Vandenberghe, and Boyd as a subroutine.  
By this we confirmed the known results for dimensions $d=2,3,4$
and extended them to $d=5$ (see Table \ref{table:latpackcov}).
Note that for the lattice covering problem there exists a local optimum 
for each triangulation, while this is not the case in higher dimensions
(see \cite{sv-2004b}) and for the lattice packing-covering problem.

Using our implementation 
we also found two lattices which currently give the best known 
covering and packing-covering in dimension $6$:

\begin{theorem}[\cite{sv-2004}]
In dimension $6$, there exits a lattice $L^{c}_6$ with 
$\Theta(L^{c}_6)= 2.4648\ldots$
and a lattice $L^{pc}_6$ with
$\gamma(L^{pc}_6)
= 1.4110\ldots$.
\end{theorem}

In \cite{sv-2004b} we show that the root lattice ${\mathsf E}_8$
does not give a locally optimal lattice covering, by constructing 
a refining triangulation $\MD$ of $\Del(Q_{{\mathsf E}_8})$ in which 
$\Theta$'s local optimum is not attained by the PQF $Q_{{\mathsf E}_8}$.
The PQF found in this way even beats the formerly best known
value $\Theta(\mathsf{A}^*_8)$ by more than $12\%$. By looking at a bistellar neighbor of the triangulation
$\MD$, we found the currently best known covering lattice in dimension $8$.

\begin{theorem}[\cite{sv-2004b}]
In dimension $8$, there exists a lattice $L^{c}_8$ with 
$\Theta(L^{c}_8)= 3.1423\ldots$.
\end{theorem}

Looking at the results in dimension $d=6,8$ it is interesting to observe that
we found the new covering lattices by looking at triangulations
refining the Delone subdivisions of the lattices ${\mathsf E}^\ast_d$. 
By looking at a corresponding refinement of ${\mathsf E}^\ast_7$,
we also found a new covering record in dimension $7$.
It remains to see if these results have a common explanation\dots

\subsection{Sufficient Conditions for Local Optima}
 
A disadvantage of finding local optima via convex programming is that
solutions can only be approximated. But this is an inherent problem:
In contrast to the lattice packing problem, local optima to the other two
problems can in general not be represented by rational numbers. One has
to use algebraic numbers instead. In some cases it might be possible,
e.g.\ with additional information on the automorphism group, to attain
exact coordinates from a first approximation (see \cite{sv-2004} for
an example).  In other cases we might have a conjectured optimal form
$Q'$ and want to compute a ``certificate'' verifying its local
optimality.  The following two propositions give such a criterion in
terms of the gradient $g_L(Q')=\gradient |\BR_L| (Q')$ of the regular
surfaces $|\BR_L(Q)|=0$ at $Q'$.  Both are a consequence of the
geometric fact that we have a local optimum at $Q'\in\SC(\MD)$ with
$\mu(Q')=1$ if and
only if there exists a hyperplane through $Q'$, separating the convex sets
$\{ Q\in\sdo : L\in \MD, |\BR_L(Q)|\geq 0\}$ and $\{Q\in\sdo :
\det(Q)\geq \det(Q')\}$, respectively $\{Q\in\sdo : \lambda(Q)\geq
\lambda(Q') \}$.

\begin{proposition}[Barnes and Dickson \cite{bd-1967}]
\label{prop:barnes-dickson}
Let $\MD$ be a Delone triangulation.  Then $Q\in\SC(\MD)$ with $\mu(Q)
= 1$ is a unique locally optimal solution to the lattice covering
problem if and only if
\[
Q^{-1} \in -\cone\{g_L(Q) : L\in \MD \mbox{ with } |\BR_L(Q)|=0 \}
.
\]
\end{proposition}

The corresponding result for the lattice packing-covering problem is

\begin{proposition}[\cite{sv-2004}] \label{prop:our-analogon}
Let $\MD$ be a Delone triangulation.  Then $Q\in\SC(\MD)$ with
$\mu(Q) = 1$ is a locally optimal solution to the lattice 
packing-covering problem if and only if
\[
\cone\{ \vec{v}\vec{v}^t : \vec{v}\in \Min(Q)\} 
\cap -\cone\{g_L(Q) : L\in \MD \mbox{ with } |\BR_L(Q)|=0 \}
\not= \emptyset
.
\]
\end{proposition}

Combining these two propositions we get

\begin{corollary}
\label{cor:combination}
Let $\MD$ be a Delone triangulation.  Then $Q\in\SC(\MD)$ is a unique
locally optimal solution to the lattice packing-covering problem, if
$Q$ is eutactic and a locally optimal solution to the lattice covering
problem.
\end{corollary}

\begin{example}
We can use Corollary \ref{cor:combination} to show that
$\mathsf{A}^*_d$, $d \geq 2$, with $\gamma(\mathsf{A}^*_d)=
\sqrt{\frac{d+2}{3}}$ is locally optimal for the lattice
packing-covering problem. Ryshkov \cite{ryshkov-1974} gave another
proof of this fact. The lattice $\mathsf{A}^*_d$ is known to give a locally optimal
lattice covering (see \cite{gameckii-1962}, \cite{gameckii-1963},
\cite{bleicher-1962}).  A quadratic form $Q_{\mathsf{A}^*_d}$
associated with $\mathsf{A}^*_d$ is
$$
Q_{\mathsf{A}^*_d} =
\begin{pmatrix}
d & -1 & \cdots & -1 \\
-1 & \ddots & \ddots & \vdots \\
\vdots & \ddots & \ddots & -1\\
-1 & \cdots & -1 & d \\
\end{pmatrix}
$$
with $\mu(Q_{\mathsf{A}^*_d})=\frac{d(d+2)}{12}$ and $\lambda(Q)=d$.  The
set $\Min(Q_{\mathsf{A}^*_d})$ contains exactly $2(d+1)$ elements,
namely the standard basis vectors $\vec{e}_1,\dots,\vec{e}_d$,
their negatives and $\pm \sum_{i=1}^d \vec{e}_i$ (see \cite{cs-1988}).
Thus in particular
$$
(d+1)\cdot Q^{-1}_{\mathsf{A}^*_d}
= Q_{\mathsf{A}_d} =
\begin{pmatrix}
2 & 1 & \cdots & 1 \\
1 & \ddots & \ddots & \vdots \\
\vdots & \ddots & \ddots & 1\\
1 & \cdots & 1 & 2 \\
\end{pmatrix}
=\sum_{\vec{v}\vec{v}^t : \vec{v}\in \Min(Q_{\mathsf{A}^*_d})}  \vec{v}\vec{v}^t
$$
Hence, $Q_{\mathsf{A}^*_d}$ is eutactic because
$Q^{-1}_{\mathsf{A}^*_d} \in \relint\cone\{ \vec{v}\vec{v}^t :
\vec{v}\in \Min(Q_{\mathsf{A}^*_d})\}$ and therefore the assertion follows.
\end{example}

Propositions \ref{prop:barnes-dickson} and \ref{prop:our-analogon}
assume that $\MD = \Del(Q)$ is a Delone triangulation.  If this is not
the case, the situation becomes more complicated, in particular for
the lattice packing-covering problem. 

For the covering problem we only have to add a condition on the set
\[
\MV_\MD = 
\bigcup_{\MD' < \MD} \{Q \in \overline{\SC(\MD')} : |\BR_L(Q)| \geq 0 \mbox{ for all } L \in \MD'\}.
\]
where $\MD' < \MD$ denotes that $\MD'$ is a Delone
triangulations refining $\MD$. This set is a subset of $\{Q \in \sdo :
\mu(Q) \leq 1\}$.  We require that $\MV_\MD$ is {\em separatable at $Q$},
that is, there exists a supporting hyperplane of $\MV_\MD$ through
$Q$.  This is in particular the case, if there exists a small $r>0$
such that $(Q+r B^{\binom{d+1}{2}})\cap \MV_\MD$ is convex.

\begin{proposition}
Let $\MD$ be a Delone subdivision and $Q\in\SC(\MD)$ with
$\mu(Q) = 1$. Then 
\begin{enumerate}
\item
$Q$ is a locally optimal solution to the lattice covering problem, if
and only if $\MV_\MD$ is separatable at $Q$ and
\[
Q^{-1} \in -\cone\{g_L(Q) : \mbox{$L\in \MD'<\MD$ with $|\BR_L(Q)|=0$} \} .
\]
\item
$Q$ is a locally optimal solution to the lattice packing-covering
problem, if $\MV_\MD$ is separatable at $Q$ and
\[
\cone\{ \vec{v}\vec{v}^t : \vec{v}\in \Min(Q)\} 
\cap 
-\cone\{g_L(Q) : \mbox{$L\in \MD'<\MD$ with $|\BR_L(Q)|=0$} \}
\neq \emptyset
.
\]
\end{enumerate}
\end{proposition}

In case of the lattice packing-covering problem the ``only if'' part
is missing, because we can not exclude the case of a locally optimal
solution $Q'$with $\MV_\MD$ not being separatable at $Q'$.  This is
due to the fact that $\{Q\in\sdo : \lambda(Q)\geq \lambda(Q') \}$ is
not smooth in contrast to $\{Q\in\sdo : \det(Q)\geq \det(Q')\}$.

This phenomenon seems to happen to PQFs associated to the root lattice
$\mathsf{E}_8$.  This lattice is known to give a globally optimal
solution to the lattice packing problem, but not a locally optimal solution to
the lattice covering problem (see \cite{sv-2004b}).  Nevertheless
computational experiments support the

\begin{conjecture}
The root lattice $\mathsf{E}_8$ gives a locally optimal solution for the 
lattice packing-covering problem.
\end{conjecture}

Zong \cite{zong-2002} even conjectured that $\mathsf{E}_8$
gives the unique globally optimal solution to the lattice
packing-covering problem in dimension $8$.

\subsection{Local Optima via Local Lower Bounds}

In \cite{sv-2004} we describe a way to attain local lower bounds for
the covering density and the packing-covering constant due to Ryshkov
and Delone.  A variant of this method is successfully used in
\cite{sv-2004b} to prove the local covering optimality of the Leech
lattice.  Here we describe a corresponding local lower bound for the
lattice packing-covering problem. As an example we use it to prove the
local packing-covering optimality of the Leech lattice directly.

\begin{proposition}
\label{prop:moments2}
Let $L_1, \ldots, L_n$ be a collection of Delone simplices of a PQF $Q$.
Then
$$\gamma(Q) \geq 2 \sqrt{\frac{\trace(F\cdot Q_{F})}{(d+1)\lambda(Q_F)}}$$
with the PQF $F=\frac{1}{n(d+1)} \sum_i\sum_{k \neq l} \vec{v}_{i,k}\vec{v}_{i,l}^t$
and a PQF $Q_F$ with
\[
F \in \cone \{\vec{v}\vec{v}^t : \vec{v}\in\Min(Q_F)\}
.
\]
\end{proposition}

One can prove this Proposition by doing obvious modifications to the
proof of Proposition~10.6 in \cite{sv-2004}. As in
Proposition~\ref{prop:our-analogon} we use the following fact: A
linear function $f(Q) = \trace(F \cdot Q)$, with a PQF $F$, has a
minimum on the homogeneous minimum $\lambda$ surface $\{Q \in \sdo :
\lambda(Q) = \lambda\}$ at $Q_F$ if and only if $F \in
\cone\{\vec{v}\vec{v}^t : \vec{v} \in \Min(Q_F)\}$.
In particular, if $Q$ is eutactic and $F = Q^{-1}$, 
then Proposition~\ref{prop:moments2} is
immediately applicable with $Q_F = Q$.

\begin{example}
We use Proposition~\ref{prop:moments2} to show that the
Leech lattice is a locally optimal packing-covering lattice.

Let us briefly review some necessary properties of the Leech lattice
$\Lambda$.  For further reading we refer to \cite{cs-1988}.  An
associated PQF $Q_{\Lambda}$ has (up to congruences) $23$ different
Delone polytopes attaining the maximum squared circumradius
$\mu(Q_\Lambda)=2$.  One of them is the Delone simplex $L$ of type
$\mathsf{A}_{24}$. 

Now we apply Proposition~\ref{prop:moments2} to the orbit of $L$ under
the automorphism group $\Co_0 = \{T \in \GL_{24}(\Z) : T^t Q_{\Lambda}
T = Q_{\Lambda}\}$ of $Q_{\Lambda}$. We get $F = \frac{1}{25|\Co_0|}
\sum_{T \in \Co_0} \sum_{\vec{e}} \vec{e}\vec{e}^t$, where $\vec{e}$
runs through all the edge vectors of $TL$. In \cite{sv-2004b} it was
shown that $F = \frac{5^2}{2^2\cdot3} Q_{\Lambda}^{-1}$.  Due to the
fact that $\Min Q_{\Lambda}$ is a {\em spherical $2$-design} with
respect to the inner product given by $Q_{\Lambda}$, we know (see
\cite{sv-2004b} for details)
\[
\sum_{\vec{v} \in \Min(Q_{\Lambda})} \vec{v}\vec{v}^t
= \frac{|\Min(Q_{\Lambda})|}{d} Q_\Lambda^{-1}
.
\]
Thus, $Q_{\Lambda}$ is eutactic and we may use $Q_F = Q_{\Lambda}$ in Proposition
\ref{prop:moments2}.  With $\lambda(Q_\Lambda)=4$ we derive
\[
\gamma(Q)
\geq
2\sqrt{\frac{\frac{5^2}{2^2\cdot3}\cdot 24}{25\cdot 4}}=\sqrt{2}=\gamma(Q_\Lambda)
\]
for all PQFs $Q$ with Delone simplices $TL$, $T \in \Co_0$,
which proofs the assertion.
\end{example}

\bibliographystyle{amsalpha}

\providecommand{\bysame}{\leavevmode\hbox to3em{\hrulefill}\thinspace}
\providecommand{\MR}{\relax\ifhmode\unskip\space\fi MR }
\providecommand{\MRhref}[2]{%
  \href{http://www.ams.org/mathscinet-getitem?mr=#1}{#2}
}
\providecommand{\href}[2]{#2}

\end{document}